\documentclass[a4paper,twoside,11pt,reqno]{amsart}

\usepackage[T1]{fontenc}
\usepackage{courier}
\usepackage[scaled=0.92]{helvet}

\usepackage{amscd,amsmath,amssymb,amsthm,amsfonts,epsfig,graphics}

\usepackage[english]{babel}
\usepackage{amsthm,amsfonts,amssymb,amsmath}
\usepackage{hyperref}
\hypersetup{
colorlinks=true,
citecolor=red,
linkbordercolor={1 1 1},
citebordercolor={1 1 1},
}

\usepackage{mathrsfs}
\usepackage{xcolor} 

\theoremstyle{plain}
\newtheorem{satz}{Theorem}[section]
\newtheorem{prop}[satz]{Proposition}
\newtheorem{cor}[satz]{Corollary}
\newtheorem{lem}[satz]{Lemma}

\theoremstyle{definition}
\newtheorem{defn}[satz]{Definition}
\newtheorem{rem}[satz]{Remark}

\newtheorem{hyp}[satz]{Hypothesis}

\newcommand{\la}{\langle}
\newcommand{\re}{\rangle}
\newcommand{\mx}{\mbox}
\newcommand{\rw}{\rightarrow}
\newcommand{\de}{\displaystyle}
\newcommand{\ml}{\mathcal}
\newcommand{\tf}{\textbf}

\newcommand{\pl}{\partial}

\newcommand{\beq}[1]{\begin{equation} \label{#1}}
\newcommand{\eeq}{\end{equation}}
\newcommand{\beqar}{\[ \begin{array}{rcl}}
\newcommand{\eeqar}{\end{array} \]}

\newcommand{\tcl}[1]{\textcolor{black}{#1}}

\providecommand{\ep}{\varepsilon}
\providecommand{\ph}{\varphi}

\providecommand{\RR}{\mathbb{R}}
\providecommand{\CC}{\mathbb{C}}

\providecommand{\ZZ}{\mathbb{Z}}
\providecommand{\NN}{\mathbb{N}}

\providecommand{\TT}{\mathbb{T}}

\newcommand{\snorm}[2]{\left| #1\right|_{#2}}
\newcommand{\norm}[2]{\left \lVert#1 \right\rVert_{#2}}

\usepackage[utf8]{inputenc}

\hypersetup{pdfauthor={Fortunati, Wiggins},%
            pdftitle={Aperiodic Nekhoroshev Theorem},%
            pdfsubject={Hamiltonian Mechanics},%
            pdfkeywords={Nekhoroshev stability, Aperiodic time dependence}%
}

\DeclareMathOperator{\dist}{dist}
\DeclareMathOperator{\spn}{span}

\DeclareMathOperator{\id}{Id}

\makeatletter
\g@addto@macro{\endabstract}{\@setabstract}
\newcommand{\authorfootnotes}{\renewcommand\thefootnote{\@fnsymbol\c@footnote}}%
\makeatother

\begin{document}
\title[Normal form and Nekhoroshev stability]{Normal Form and Nekhoroshev stability for nearly-integrable Hamiltonian systems with unconditionally slow aperiodic time dependence}
\author{Alessandro Fortunati}
\thanks{This research was supported by ONR Grant No.~N00014-01-1-0769 and MINECO: ICMAT Severo Ochoa project SEV-2011-0087.}
\address{School of Mathematics, University of Bristol, Bristol BS8 1TW, United Kingdom}
\email{alessandro.fortunati@bristol.ac.uk}
\keywords{Hamiltonian systems, Nekhoroshev theorem, Aperiodic time dependence.}
\subjclass[2000]{Primary: 70H08. Secondary: 37J25, 37J40}

\author{Stephen Wiggins}
\email{s.wiggins@bristol.ac.uk}

\maketitle

\begin{abstract}
The aim of this paper is to extend the results of Giorgilli and Zehnder for aperiodic time dependent systems to a case of general nearly-integrable convex analytic Hamiltonians. The existence of a normal form and then a stability result are shown in the case of a slow aperiodic time dependence that, under some smallness conditions, is independent of the size of the perturbation. 

\end{abstract}
\section{Introduction}
The study of the solutions of a near-integrable Hamiltonian system goes back to Poincaré \cite{poi}, who emphasized the relevance of this model by describing it as ``General problem of the Dynamics''. Motivated by problems arising from Celestial Mechanics, stability would have been among the most interesting (and urgent) questions to be addressed.\\
The persistence of invariant tori of an integrable Hamiltonian under small perturbations, initially faced by Kolomogorov \cite{kolm}, gave a powerful answer to this problem: the \emph{perpetual stability} of certain invariant tori. A different proof, due to Arnold  \cite{arn1}, showed that the invariant tori persisted 
on a very special subdomain of the phase space $\ml{D}$,  a Cantor set whose (Lebesgue)  measure is ``close'' to the (Lebesgue) measure of the phase space   $\ml{D}$.  As a drawback of Arnold's method of proof, the construction of an ``exact'' normal form (i.e. by an infinite number of steps) provided by this result, holds in a very special subdomain of the phase space $\ml{D}$ (Cantor set), which measure is close to the $\ml{D}$ one, but it is completely different from a topological point. This ``high probability'' \cite{poes2} to lie on an invariant torus is not adequate for certain applications.\\ 
The possibility of a weaker statement on a more suitable domain from the applications point of view, consisted of the so-called \emph{effective stability}. After the initial contributions by Moser \cite{moser55}, Littlewood \cite{littlewood59a}, \cite{littlewood59b} and subsequently by Glimm \cite{gli}, it was realized in a general setting by Nekhoroshev \cite{nekhoroshev77}. The starting point changed the KAM point of view: by keeping only a order $r-$truncated (resonant) normal form, it is possible to preserve an open subset of the phase space, then a careful choice of $r$ can be made in order to obtain a stability time as large as possible. Obviously this result, as it is, is only of local nature. The decisive contribution of Nekhoroshev was the so-called \emph{geometric part}, in which the entire phase space is covered by using suitable subsets having known resonance properties (geography of resonances) to which the normal form result can be applied.\\
The relevance of this result has rapidly raised the interest of the scientific community outside the Russian school, especially in Italy with  Benettin, Galgani, Gallavotti and Giorgilli, e.g. \cite{bgg85}, \cite{galla} and subsequent papers, then in France with Lochak, e.g. \cite{lo}, who developed a new approach (simultaneous Diophantine approximations) able to enlarge the exponent of the stability bound. The steepness feature of the unperturbed Hamiltonian initially considered by Nekhoroshev, was profitably replaced by the (slightly) less general but remarkable \emph{convexity} in the '80, and then weakened to \emph{quasi-convexity} in \cite{poes}.
As it was reasonable to expect, solutions ``close'' to a KAM torus would  possess special stability properties. This aspect was made precise in  \cite{morbidelli} and \cite{perry}, with the former reference  showing that  solutions starting  exponentially close to a KAM torus are indeed \emph{super-exponentially} stable. \\
Meanwhile, the paper \cite{gz92} proposed a different direction, by considering a model of the form $H(x,y,t)=|\dot{y}|^2/2+f(x,y,t)$ (i.e. convex). The dependence of $f$ is quasi-periodic on $x$ but only analytic on $t$, introducing in this way, for the first time, an aperiodic time dependence. A Nekhoroshev type result is shown for motions with ``high'' kinetic energy, i.e., after a time rescaling $t=:\ep \tau$, for small $\ep$ and bounded energies. As a side effect, the dependence on $\tau$ turns out to be \emph{slow} with $\ep$.\\
The key property used in the perturbative setting of \cite{gz92}, consists of the possibility to disregard the dependence on $\tau$ in the solution of the homological equation. This has the irrelevant effect of losing t control of the variable canonically conjugate to $\tau$, say $\eta$, the latter being  a fictitious variable\footnote{Furthermore if $t$ is a scalar variable as in our case, bounded variations of the actions imply that $\eta$ is bounded \tcl{as well}, simply by the conservation of energy.}. This argument of \emph{partial normal form} substantially simplifies the discussion and, as it will be shown, allows an immediate interfacing with the standard quasi-periodic case. \\
Despite these innovative features, the mentioned result has not been widely known for more than twenty years. Recently, the problem has been reconsidered in \cite{boun13}, giving an outline of the elements necessary to adapt previous results by the same author to a system of the form $H(I,\theta,t)=h(I)+\ep f(I,\theta,\ep^c t)$, $c \in \RR^+$. A slow time dependence similar to \cite{gz92} is considered. \\
Also in the light of the  applications of this kind of result pointed out in \cite{wiggins}, the aim of this paper is to extend the results by \cite{gz92} to  more general systems, along the lines of a proof of the Nekhoroshev theorem described in the comprehensive paper \cite{giorgilli02}. The fully constructive setting given by the Lie transform method allowed a deeper analysis of the slow time dependence problem. As also stressed in \cite{boun13}, an hypothesis of slow time dependence is completely reasonable. Otherwise, it is natural to expect the existence of \emph{ad hoc} perturbations able to ``drive'' the solutions along some resonance channel. Roughly, the role of the small parameter is to create a safe separation between the frequencies of the unperturbed system and those produced by the perturbation. Nevertheless, by considering a two-parameters system (\ref{eq:ham}), we show that the ``speed'' of the time dependence and the size of the perturbation, under some smallness hypothesis, should  not be necessarily related. This is the feature behind the \emph{unconditionally slow dependence}. It leads to a great advantage from the applications point of view and is an extension of the results obtained in \cite{gz92} and \cite{boun13}.\\
The partial normal form, whose existence is shown in sec. $3$, allows to use \emph{exactly} the same geometric arguments of a standard Nekhoroshev result. For this reason, more general hypothesis on the unperturbed Hamiltonian (used only in the geometric part) than the convexity are not addressed here and a brief outline of the argument described in \cite{giorgilli02}, is given in sec. $4$ for the sake of completeness.

\section{Set-up and main result}
Let $\ml{G}$ be an open subset of $\RR^n$ and consider the nearly integrable system described by the following Hamiltonian function 
\beq{eq:ham}
H(I,\ph,t):=h(I)+\ep f(I,\ph,\mu t) \mx{,}
\eeq
where $I=(I_1,\ldots,I_n) \in \mathcal{G}$, $\ph=(\ph_1,\ldots,\ph_n) \in \TT^n$ is a set of 
action-angle variables and $t \in \RR$ is an additional variable (time) on which $f$ \tcl{does not need to depend quasi-periodically}.\\ 
As usual, by setting $\xi:=\mu t$ and denoting by $\eta \in \RR$ the variable conjugate to $t$, Hamiltonian 
(\ref{eq:ham}) can be seen as autonomous in the extended phase space 
$\mathcal{D}:= \mathcal{G} \times \RR \times \TT^n \times \RR \ni (I,\eta,\ph,\xi) $ in the form
\beq{eq:hamnotempo}
H(I,\ph,\xi,\eta):=h(I)+\mu \eta+\ep f(I,\ph,\xi) \mx{.}
\eeq
Given two control parameters $\rho,\sigma \in (0,1]$, consider the complex neighbourhood of $\ml{D}$, defined as
$\ml{D}_{\rho,2 \sigma}:=\ml{G}_{\rho} \times 
\ml{R}_{\rho} \times \TT_{2 \sigma}^n \times R_{2 \sigma}$ where
\beqar
\mathcal{G}_{\rho}&:=& \de \bigcup_{I \in \ml{G}} \Delta_{\rho}(I),\qquad 
\Delta_{\rho}(I):=\{\hat{I} \in \CC^n : |\hat{I}-I|<\rho\}\\
\ml{R}_{\rho}&:=&\{\hat{\eta} \in \CC : |\Im \hat{\eta}|<\rho\}\\
\TT_{2 \sigma}^n &:=& \{\hat{\ph} \in \CC^n: |\Im \hat{\ph}| < 2 \sigma\}\\
R_{2 \sigma}&:=& \{\hat{\xi} \in \CC: |\Im \hat{\xi}| <  2 \sigma\}
\eeqar  
The space $\mathcal{D}_{\rho,2 \sigma}$ is endowed with the usual supremum norm 
\[
\snorm{F}{\rho,\sigma} := \sup_{z \in \mathcal{D}_{\rho,2 \sigma}} |F(z)| \mx{.}
\] 
For any analytic functions $F=F(I,\ph,\xi) \in \ml{D}_{\rho,2 \sigma}$, admitting a Fourier expansion of the form 
\[
 F(I,\eta,\xi)=\sum_{k \in \ZZ^n} f_k(I,\xi) e^{i k \cdot \ph} \mx{,}
\]
the \emph{Fourier norm} is defined as
\[
\norm{F}{\rho,\sigma}:=\sum_{k \in \ZZ^n} 
\snorm{f_k}{\rho,\sigma}e^{|k| \sigma} \mx{,}
\]
where $|k|=|k_1|+\ldots+|k_n|$. System (\ref{eq:ham}) will be studied under the following 
\begin{hyp}\label{hyp:reg} $h(I)$ and $f(I,\ph,\xi)$ are holomorphic and bounded functions on the
space $\mathcal{D}_{\rho,2 \sigma}$, in particular 
\beq{eq:flimitata}
 \snorm{f}{\rho,2 \sigma} =:C_f < +\infty \mx{.}
\eeq
\end{hyp}
\begin{hyp}\label{hyp:regdue}
The unperturbed Hamiltonian $h(I)$ is a convex function, i.e. there exists two constants $M \geq m >0$ such that, for all $I \in \ml{G}_{\rho}$
\beq{eq:nondegeneracy}
|\pl_I^2 h(I) v|\leq M |v|, \qquad |\la \pl_I^2 h(I) v,v \re| \geq m |v|^2 \mx{,}
\eeq
for all $v \in \RR^n$.
\end{hyp}
The set of parameters $\rho,\sigma,M,m,C_f$ are characterized by  a given Hamiltonian and will be supposed fixed once and for all. Let us define
\beq{eq:ftildeeco}
\tilde{\ml{F}}=C_f \left( 
\frac{1+e^{-\frac{\sigma}{2}}}{1-e^{-\frac{\sigma}{2}}}\right)^n,\qquad \lambda_{\ep,\mu}:=\mu+e \tilde{\ml{F}} \ep \mx{.}
\eeq
In this framework, the main result is stated as follows
\begin{satz}[Aperiodic Nekhoroshev]\label{thm:aper}
Assume hypotheses \ref{hyp:reg} and \ref{hyp:regdue}. Then there exists constants $\Delta^*$ and $\ml{T}$ depending on
$\rho,\sigma,M,m,C_f$ and $n$ such that, if $\ep$ and $\mu$ satisfy
\beq{eq:condizioneepsilon}
\lambda_{\ep,\mu}<\frac{1}{3^4 \Delta^*} \mx{,}
\eeq
then orbits $(I(t),\ph(t))$ of (\ref{eq:ham}) starting in $\ml{G} \times \TT^n$ at $t_0$, satisfy 
\[
 |I(t)-I(t_0)|<(\Delta^* \lambda_{\mu,\ep})^{\frac{1}{4}} \rho, \qquad \mx{for} \qquad
 |t-t_0|<\frac{\ml{T}}{\ep}\exp \left[ \left( \frac{1}{\Delta^* \lambda_{\mu,\ep}}\right)^{\frac{1}{2a}}\right] \mx{,}
\]
where $a=n^2+n$.
\end{satz}
We remark that $\Delta^*$ is defined in sec. $4$ .
The main feature of this formulation is that the smallness condition (\ref{eq:condizioneepsilon}) allows a certain freedom in the choice of $\ep$ and $\mu$. Essentially, within the described threshold, parameters $\ep$ and $\mu$ can be treated as \emph{independent}. In principle, $\mu$ is even allowed to be increased, as $\ep$ tends to $0$, still preserving a normal form result (see Theorem \ref{thm:normal}) and a stability estimate, despite the fact that, in such case, it is a worse estimate as the bound is only $O(\ep^{-1})$ as $\ep$ vanishes. In any case, this is an extension of the already existing results, e.g. \cite{gz92} and \cite{boun13}, in which $\mu=\ep^c$, $c>0$, i.e. the aperiodic time dependence is forced to be slower as $\ep$ gets smaller. Consistently, Theorem \ref{thm:aper} coincides with the result stated in \cite{giorgilli02}, when the limit $\mu \rw 0$ is considered.\\
As mentioned above, the proof of this theorem is deduced from \cite{giorgilli02} and \cite{gz92} with some technical modifications that will be discussed in detail. The same notation of \cite{giorgilli02} is preserved as much as possible, for a more efficient comparison between the results. 

\section{Normal form}
Throughout this section, $\rho$ will be replaced with $\delta$ in order to avoid confusion in the final estimate. 
\subsection{Basic notions and statement}
\begin{defn}
A subset $\ml{M}$ of $\ZZ^n$ is said to be a \tf{resonance module} if it satisfies
\beq{eq:condmod}
\spn(\ml{M}) \cap \ZZ^d = \ml{M} \mx{.}
\eeq
\end{defn}
If $(k_1,\ldots,k_n) \in \ZZ^n$ is a basis for $\ml{M}$ (i.e. $\ml{M}=\{\alpha_1 k_1+\ldots+\alpha_n k_n:\alpha_j \in
\ZZ\}$), with $\spn(\ml{M})$ we denote the set $\{\gamma_1 k_1+\ldots+
\gamma_n k_n: \gamma_j \in \RR\} \subseteq \RR^n$. 
Hence, the purpose of condition (\ref{eq:condmod}), is to exclude subspaces of $\ZZ^d$ which contain less points of the lattice than the real space $\spn(\ml{M})$.
\begin{defn}Let $\ml{M}$ be a resonance module, $\alpha \in \RR^+$ and $N \in \NN$. A
subspace $\ml{V}$ of the action space $\ml{G}$ is said to be a \tf{non-resonance domain} of the type
$(\ml{M},\alpha,\delta,N)$ if, for all $k 
\in \ZZ^n \setminus \ml{M}$ with $|k|<N$, the following condition holds 
\[
 |\langle \omega(I),k \rangle |> \alpha,\qquad \forall I \in \ml{G}_{\delta} \mx{,}
\]
where $\omega(I):=\pl_I h(I)$. 
\end{defn}
It will be denoted by $\tilde{\ml{D}}_{\delta,\sigma} \subseteq \ml{D}_{\delta,\sigma}$, the complex extension of $\ml{D}$ with $\ml{G}$ is replaced with $\ml{V}$.
\begin{satz}[Existence of a normal form]\label{thm:normal}  
Consider the Hamiltonian (\ref{eq:hamnotempo}) with the regularity assumptions of hypothesis \ref{hyp:reg} and the associated parameters $\delta,\sigma,C_f$. Given $\ml{M}$ a resonance module, $r,K \in \NN$ and $\alpha \in \RR^{+}$, suppose that $\ml{V} \subseteq \ml{G}$ is a non-resonance domain of the type $(\ml{M},\alpha,\delta,N)$ with $N=rK$,  and that $r,K,\ep$ and $\mu$ are such that
\begin{subequations} 
\begin{align}
\Delta:=\frac{2^{8} r }{\alpha \delta \sigma} \lambda_{\ep,\mu} +4 e^{-K \frac{\sigma}{2}} & \leq \frac{1}{2}
\mx{,}
\label{eq:smallness}\\
e^{-K \frac{\sigma}{2}} & \geq (8+3e^2)^{-1} \mx{.}\label{eq:kappa}
\end{align} 
\end{subequations}
Then there exists a symplectic, $\ep-$close to identity, analytic change of variables $z \rw \ml{C}_r(z)$ defined
on
$\tilde{\ml{D}}_{\frac{3}{4}(\delta,\sigma)}$ such that 
\[
\tilde{\ml{D}}_{\frac{5}{8}(\delta,\sigma)}  \subset \ml{C}_r \tilde{\ml{D}}_{\frac{3}{4}(\delta,\sigma)} \subset 
\tilde{\ml{D}}_{\frac{7}{8}(\delta,\sigma)} \mx{,}
\]
(the same for holds for $\ml{C}_r^{-1}$), and casting the Hamiltonian (\ref{eq:hamnotempo}) in \tf{resonant normal form} up to order $r$, i.e.
\beq{eq:normformth}
H( \ml{C}_r (z))=h(I)+\eta+Z^{(r)}+\ml{R}^{(r+1)} \mx{,}
\eeq
with $Z^{(r)}=\sum_{k \in \ml{M}} z_k(I,\xi)e^{i k \cdot \ph}$ for all $|k|\leq N$  and
\beq{eq:stimaresto}
\norm{\ml{R}^{r+1}}{\frac{3}{4}(\delta,\sigma)} \leq 8\ep \tilde{\ml{F}} \Delta^r \mx{.}
\eeq 
\end{satz}
\begin{rem} As mentioned before, by condition (\ref{eq:smallness}) it is evident that the normal form exists as long as $\lambda_{\ep,\mu}$ is sufficiently small, no matter if there is relation or not between $\ep$ and $\mu$. The technical hypothesis (\ref{eq:kappa}) does not appear in \cite{giorgilli02}, anyway, it will be shown that it can be assumed without loss of generality. 
\end{rem}
The proof of the already stated result goes along the lines of \cite{gz92} and \cite{giorgilli02}, and can can be achieved in two step. In the first one, a suitable perturbative algorithm is
built up in order to remove the effect of the perturbation up to a pre-fixed order $r$ (except for a particular set of
harmonics given by $\ml{M}$). The perturbative scheme is initially discussed at a formal level
(i.e. disregarding the problem of the series convergence) and it is based on the Lie transform method. The subsequent step consists of giving quantitative estimate on the convergence of the scheme by using standard analytic tools (see \cite{giorgilli02}). 
\subsection{Perturbative setting: formal scheme}
Given $K$, define, for all $j \geq 1$, the following class of functions on $\tilde{\ml{D}}_{\delta,\sigma}$ 
\[
\ml{P}_j:=\{g:g(I,\ph,\xi)=\sum_{|k|<j K} g_k(I,\xi) e^{i k \cdot \ph}\} \mx{.}
\]
By setting, for all $s=1,2,\ldots,$ 
\[
H_s:= \ep \sum_{ (s-1) K \leq |k| < s K} f_k(I,\xi) e^{i k \cdot \ph} \mx{,}
\]
the Fourier expansion of $f$ had been split in such a way $H_s \in \ml{P}_s$ for all $s \geq 1$, and the Hamiltonian
(\ref{eq:hamnotempo}) reads as
\beq{eq:hamsomma}
H(I,\ph,\xi,\eta)=h(I)+\eta+ H_1+H_2+\ldots \mx{.}
\eeq
The aim is to find a local, $\ep-$close to the identity, symplectic diffeomorphism 
casting the Hamiltonian (\ref{eq:hamnotempo}) into the form (\ref{eq:normformth}). 
This is achieved via a suitable choice of a (finite) sequence of 
functions $\chi^{(r)}:=\{\chi_s\}_{s=1,\ldots,r}$ (generating sequence), and setting $\ml{C}_r \equiv T_{\chi^{(r)}}$, where $T_{\chi^{(r)}}$ is the \emph{Lie
transform} operator associated to $\chi^{(r)}$
\beq{eq:lietransform}
T_{\chi^{(r)}}:=\de \sum_{s=1}^r E_s,\qquad 
E_s:=
\left\{
\begin{array}{lcl}
\id & \quad& s=0\\
\de \frac{1}{s} \sum_{j=1}^s j \ml{L}_{\chi_s} E_{s-j} & \quad & s \geq 1
\end{array}
\right.
\eeq
and $\ml{L}_{\chi_s} g:=\{g,\chi_s\}$ stands for the \emph{Lie derivative}. \tcl{Note that if $f,g$ are two functions independent of $\eta$ (as the objects involved in the above perturbative scheme), the parenthesis $\{f,g\}$ reduces to $\sum_{i=1}^n(\pl_{\ph_i}f \pl_{I_i}g - \pl_{I_i}f \pl_{\ph_i}g )$. In particular, $\{\xi,f\}=0$ for all $f=f(I,\ph,\xi)$, i.e. $T_{\chi^{(r)}}(\xi)=\xi$. Hence the considered transformation does not act on time}.\\ 
Taking into account of (\ref{eq:hamsomma}) and writing $Z^{(r)}=Z_1+\ldots+Z_r$ and  one gets the following 
\begin{prop}
Equation (\ref{eq:normformth}) is equivalent to the following hierarchy of homological equations
\beq{eq:homologicaltrue}
\ml{L}_h \chi_s + Z_s = \psi_s \mx{,}
\eeq
for $s=1,\ldots,r$ with $Z_s \in \ml{P}_s$ and 
\beq{eq:psis}
\psi_s:=\left\{
 \begin{array}{lcl}
 H_1 &\qquad& s=1\\
 H_s+ \mu E_{s-1} \eta + \de \frac{1}{s}\sum_{j=1}^{s-1} j[\ml{L}_{\chi_j}H_{s-j}+E_{s-j}H_j] &\qquad& 2 \leq s \leq r
 \end{array}
 \right. 
\eeq
\end{prop}
\proof(Sketch). Use the well known identity
$H(T_{\chi^{(r)}})=T_{\chi^{(r)}}H$, then substitute into (\ref{eq:normformth})
the involved objects in form of (finite) sums, equating, at the $s-$th stage, terms on the same level $\ml{P}_s$. See \cite{giorgilli02} for the details.
\endproof
\begin{rem}
Note that the operator $\ml{L}_{h}$ is exactly the same as for the standard (quasi-periodic) case. This is the main advantage in considering a \textbf{partial normal form} and the solution of the homological equation (\ref{eq:homologicaltrue}) can be done by a standard comparison of Fourier coefficient. Note that in this case, at each stage, the averaged term does not depend only on $I$ but also on $\xi$. The key fact used in \cite{gz92} is that one only needs a partial normal form and this term can be anyway included in $Z$ as this does not affect the evolution of the variables $I$, but only of $\eta$. \tcl{Consistently, if the aperiodic dependence on $\xi$ is supposed to be quasi-periodic, the dependence on the ``angle'' $\xi$ is annihilated by averaging and the partial normal form becomes ``full''}.\\
Hence, the argument is reduced to the control of the extra-term $E_{s-1} \eta$ arising from the aperiodic time dependence.   
\end{rem}
As for the solution of (\ref{eq:homologicaltrue}), denote by $z_k^{(s)}$, $c_k^{(s)}$ and $v_k^{(s)}$ the Fourier coefficients of $Z_s$, $\chi_s$ and $\psi_s$ respectively. As the resonance module $\ml{M}$ is fixed, if $k=0$ or $k \in \ml{M}$ one sets $c_k^{(s)}=0$ and $z_k^{(s)}=v_k^{(s)}$. Otherwise, if only $I \in \ml{V}$ are considered, the quantity $\la k,\omega(I)\re$ is bounded away from zero and it is possible to set $c_k^{(s)}=i(\la k,\omega(I)\re)^{-1}v_k^{(s)}$ then $z_k^{(s)}=0$. This yields immediately the following two inequalities
\beq{eq:inequalities}
\norm{Z_s}{(1-d)(\delta,\sigma)} \leq \norm{\psi_s}{(1-d)(\delta,\sigma)},\qquad 
\norm{\chi_s}{(1-d)(\delta,\sigma)} \leq \frac{1}{\alpha}\norm{\psi_s}{(1-d)(\delta,\sigma)} \mx{,}
\eeq 
valid for all $d \in (0,1)$.

\subsection{Convergence}
\begin{lem}\label{lem:tecnico}
Assume hypothesis \ref{hyp:reg}, then the following sequence of ``nested'' statements holds:
\begin{enumerate}
 \item There exists $h>0$ and $\ml{F} \geq 0$ such that
 \beq{eq:hs}
   \norm{H_s}{(\delta,\sigma)} \leq h^{s-1} \ml{F},\qquad s \geq 1\mx{.}
 \eeq
 \item Supposing (\ref{eq:hs}), holds 
  \beq{eq:psiseq}
 \norm{\psi_s}{(1-d)(\delta,\sigma)} \leq \ml{F} \frac{b^{s-1}}{s},\qquad s \geq 1\mx{,}
 \eeq
 for all $d<1/4$ and some $b \geq 0$. Hence, by (\ref{eq:inequalities}), the truncated series $\sum_{j=1}^r \chi_s $ and $\sum_{j=1}^r Z_s $ are well defined on $\ml{V}$, yielding respectively $\chi^{(r)}$ and $Z^{(r)}$ as a solution of (\ref{eq:normformth}). 
 \item Assume (\ref{eq:psiseq}) and that for all $d \in (0,\frac{1}{4})$ the condition
 \beq{eq:minunmezzo}
  \frac{2 e \ml{F}}{ d^2 \alpha \delta \sigma}+b \leq \frac{1}{2} \mx{,}
 \eeq
is satisfied, then the operator $T_{\chi^{(r)}}$ (and its inverse $T_{\chi^{(r)}}^{-1}$) define a canonical
transformation on the domain $\tilde{\ml{D}}_{(1-d)(\delta,\sigma)}$ with the following properties
\[
\begin{array}{c}
\tilde{\ml{D}}_{(1-2d)(\delta,\sigma)}  \subset T_{\chi^{(r)}} \tilde{\ml{D}}_{(1-d)(\delta,\sigma)} \subset 
\tilde{\ml{D}}_{(\delta,\sigma)}\\
\tilde{\ml{D}}_{(1-2d)(\delta,\sigma)}  \subset  T_{\chi^{(r)}}^{-1} \tilde{\ml{D}}_{(1-d)(\delta,\sigma)} \subset 
\tilde{\ml{D}}_{(\delta,\sigma)}
\end{array} \mx{.} 
\]
\end{enumerate}
\end{lem}
The proof of $1$ can be found in \cite[Pag. 851]{gz92} where in addition it is shown that
\beq{eq:feh}
\ml{F}=\ep \tilde{\ml{F}}, \qquad h=e^{-K \frac{\sigma}{2}} \mx{,}
\eeq
while the third one is straightforward from a general result on the convergence of Lie transform, see \cite{giorgilli02}. The statement $2$ requires a further analysis with respect to the existing case, due to the presence of the extra-term $E_{s-1} \eta$ related to the aperiodic time dependence.\\ 
In order give an estimate for $\norm{\psi_s}{(1-d_s)(\delta,\sigma)}$, Lie operators appearing in (\ref{eq:psis}) can be treated via standard Cauchy tools, whereas a domain restriction is provided. \tcl{For, let $\tilde{\ml{D}}_{(1-d_s)(\delta,\sigma)}$ be the domains sequence}, where a convenient choice is given by
\beq{eq:ds}
 d_s:=d \sqrt{\frac{s-1}{r-1}} \mx{,}
\eeq
for all $s=1,\ldots,r$ and where $d$, appearing in (\ref{eq:minunmezzo}), will be determined later.\\
In order to control terms appearing in (\ref{eq:psis}) for all $s$, three sequences are considered that are implicitly defined by the following inequalities
\begin{subequations} 
\begin{align}
\norm{\psi_s}{(1-d_s)(\delta,\sigma)} &\leq \beta_s \ml{F} \label{stimea}\\
\norm{E_s H_j}{(1-d_{s+j})(\delta,\sigma)} & \leq  \tilde{\theta}_{s,j} \ml{F} 
\qquad j=1,\ldots,r-s \label{stimeb}
\\
\norm{E_{s-1} \eta}{(1-d_s)(\delta,\sigma)} & \leq  \gamma_{s-1} \ml{F} \label{stimec}
\end{align} 
\end{subequations}
Taking into account of (\ref{eq:hs}) and the definition of $\psi_1$, inequality (\ref{stimea}) computed for $s=0$ immediately gives $\beta_1=1$. On the other hand, recalling that $E_0=\id$, by (\ref{stimeb}) one gets $\tilde{\theta}_{0,j}:=h^{j-1}$. As for (\ref{stimec}), a Cauchy estimate, the second of (\ref{eq:inequalities}) and then (\ref{stimea}), yield
\beq{eq:stimaeuno}
\norm{E_1 \eta}{(1-d_2)(\delta,\sigma)}  \leq \norm{\frac{\pl \chi_1}{\pl \xi}}{(\delta,\sigma)} \leq \frac{1}{d \sigma} \norm{\chi_1}{(\delta,\sigma)} \leq \frac{\ml{F}}{\alpha d \sigma} \mx{,}
\eeq
i.e., defining
\beq{eq:gammagrande}
\Gamma:= (\alpha d \sigma)^{-1} \mx{,}
\eeq
one obtains $\gamma_1 =\Gamma$. Setting $m=s-1$, it is straightforward by definition that for all $m \geq 1$
\beq{eq:bydef}
\norm{E_{m} \eta}{(1-d_{m+1})(\delta,\rho)}  \leq \frac{1}{m}
\sum_{l=1}^{m-1} l \norm{\mathcal{L}_{\chi_l} E_{m-l} \eta}{(1-d_{m+1})(\delta,\rho)}+\norm{\ml{L}_{\chi_m} \eta }{(1-d_{m+1})(\delta,\sigma)} \mx{.}
\eeq
The first term of the right hand side can be estimated by using the following result (see \cite[Pag. 853]{gz92} for the proof)  
\begin{lem}\label{lem:two}
 Let $d',d'' \in \RR^+$ such that $d'+d''<1$ and 
 \begin{enumerate}
  \item $G(I,\ph,\xi)$ be analytic and bounded in $\tilde{\ml{D}}_{(1-d')(\delta,\sigma)}$,
  \item $F(I,\ph,\xi)$ be analytic and bounded in $\tilde{\ml{D}}_{(1-d'')(\delta,\sigma)}$.
 \end{enumerate}
 Then, for all $0<d<1-d'-d''$ the following inequality holds
 \beq{eq:twoparameter}
 \norm{\ml{L}_{G} F}{(1-d-d'-d'')(\delta,\sigma)} \leq C \norm{F}{(1-d'')(\delta,\sigma)} 
\norm{G}{(1-d')(\delta,\sigma)}\mx{,}
 \eeq
 where $C=2 [e (d+d')(d+d'') \delta \sigma]^{-1}$.
\end{lem}
By construction $\chi_l$ is analytic on $\tilde{\ml{D}}_{(1-d_l)(\delta,\sigma)}$ 
while $E_{m-l} \eta$ is analytic on $\tilde{\ml{D}}_{(1-d_{m-l})(\delta,\sigma)}$ 
hence by lemma (\ref{lem:two}) with $d'=d_l$ and $d''=d_{m-l}$, one gets on\footnote{Recalling (\ref{eq:ds}), one checks immediately that $ d_l+d_{m-l}< 1-d$ as required by Lemma \ref{lem:two} if $d \leq 1/3$, justifying in this way the assumption on $d$ in Lemma \ref{lem:tecnico}. Further estimates (\ref{eq:stimeuno}) and (\ref{eq:stimedue}) hold under the same assumption.}$\tilde{\ml{D}}_{(1-d_{m+1})(\delta,\sigma)}$
\beq{eq:esdue}
\norm{\mathcal{L}_{\chi_l} E_{m-l} \eta}{(1-d_{m+1})(\delta,\sigma)}
 \leq  \de \frac{2}{e  \delta (d_{m+1}-d_l)(d_{m+1}-d_{m-l})} \frac{\beta_l \ml{F}^2 }{\alpha} 
\gamma_{m-l} \mx{,}
\eeq
having used the second of (\ref{eq:inequalities}) then (\ref{stimea}) and (\ref{stimec}).\\
In conclusion, the first term of the r.h.s. of (\ref{eq:bydef}), can be bounded by (\ref{eq:esdue}), recalling (\ref{eq:ds}) and the elementary inequality
\[
(\sqrt{m}-\sqrt{l-1})(\sqrt{m}-\sqrt{m-l-1}) \geq \frac{1}{2},\qquad 1 \leq l \leq m-1 \mx{.}
\]
As for the second term, the same argument of (\ref{eq:stimaeuno}) can be used, yielding, for all $m \geq 1$
\[
\norm{E_{m} \eta}{(1-d_{m+1})(\delta,\rho)}  \leq 
\frac{C_r}{m} 
\sum_{l=1}^{m-1}  
l \beta_{l} \gamma_{m-l}+ \beta_m \Gamma \ml{F} \mx{,}
\]
where
\beq{eq:cr}
C_r:=\frac{4(r-1) \ml{F}}{ \alpha e d^2 \delta \sigma} \mx{.}
\eeq
This estimate can be clearly written in the form (\ref{stimec}) provided
\beq{eq:gammasmuno}
\gamma_{s-1}:=\frac{C_r}{s-1} \sum_{l=1}^{s-2} l \beta_l \gamma_{s-l-1} 
+ \Gamma \beta_{s-1}
\eeq
The remaining estimates (see \cite{giorgilli02} for the details), can be found in the same way and take the form
\begin{subequations}
\begin{align}
\norm{\mathcal{L}_{\chi_l} H_{s-l}}{(1-d_s)(\delta,\sigma)}  & 
\leq  \de C_r \beta_l h^{s-l-1} \ml{F} \label{eq:stimeuno}\\
\norm{E_s H_j}{(1-d_{s+j})(\delta,\sigma)} & \leq 
\de \frac{C_r}{s} \sum_{l=1}^s l \beta_l \tilde{\theta}_{s-l,j} \ml{F} \label{eq:stimedue}
\end{align}
\end{subequations}
Collecting the estimates, one obtains the following system of recurrence equations
\beq{eq:seqdue}
\left\{
\begin{array}{rcl}
\beta_s&:=& \de  h^{s-1}+\mu \gamma_{s-1}+\frac{1}{s}  
\sum_{l=1}^{s-1} l 
h^{l-1} \theta_{s-l}+\frac{C_r}{s} \sum_{l=1}^{s-1} l \beta_l h^{s-l-1}\\
\theta_{s}&:=&\de \frac{C_r}{s}\sum_{l=1}^s l \beta_l \theta_{s-l}\\
\gamma_{s} & = & \de \frac{C_r}{s} \sum_{l=1}^{s-1} l \beta_l \gamma_{s-l}+ \Gamma \beta_s  \\
\end{array}
\right.
\eeq
(see below for the definition of $\theta_s$) with the following initial conditions  
\[
\beta_1=1,\qquad \theta_0=1,\qquad \gamma_1=\Gamma \mx{.}
\]
Indeed, a comparison between (\ref{eq:stimedue}) and (\ref{stimeb}) gives $\tilde{\theta}_{s,j}=(C_r/s)\sum_{j=1}^s j \beta_j \tilde{\theta}_{s-j,j}$. From the latter, by using $\tilde{\theta}_{0,j}=h^{j-1}$ and defining $\theta_s:=\tilde{\theta}_{s,1}$, it is easy to check that $\tilde{\theta}_{s,j}=h^{j-1} \theta_s$. \tcl{This leads to the second equation. As for the first one, starting from (\ref{eq:psis}), use hypothesis (\ref{eq:hs}), (\ref{stimec}) compared with (\ref{eq:gammasmuno}), (\ref{stimeb}) and again the already obtained expression for $\tilde{\theta}_{s,j}$}. \\ 
\\
The following result provides the control of the behaviour of $\beta_s$.  
\begin{lem}\label{lem:minun} Suppose $2C_r \leq 3 h$. Then condition (\ref{eq:psiseq}) holds provided  
\beq{eq:minoreuno}
b:=4(h+C_r+\mu \Gamma) < 1\mx{.}
\eeq
\end{lem}
\begin{rem}
If $\mu=0$ this choice coincides with $b$ found in \cite{giorgilli02}. 
\end{rem}
\proof
Condition (\ref{eq:psiseq}) is proven if 
\beq{eq:betacond}
\beta_s \leq \frac{\tau^{s-1}}{s} \mx{,}
\eeq 
for all $s$ and some $\tau  \leq b $. Note that this is trivially true for $s=1$, hence suppose it for all $s \leq r-1$ and proceed by induction. 
\begin{prop}\label{prop:dueotto} Suppose (\ref{eq:betacond}), then 
\beq{eq:seqces}
\theta_s  \leq  \de \frac{C_r}{s} (\tau+ C_r)^{s-1}, \qquad 
\gamma_s  \leq  \de \frac{\Gamma}{s} (\tau + C_r)^{s-1} \mx{.}
\eeq
\end{prop}
\proof Define $\hat{\theta}_s := C_r \sum_{j=1}^s \tau^{j-1} \hat{\theta}_{s-j}$ with $\hat{\theta}_0=1$. By (\ref{eq:betacond}), follows from (\ref{eq:seqdue}) that $\theta_s \leq \hat{\theta}_s/s$ for all $s$. By a slight variant of the argument used for Catalan numbers (see e.g. \cite{kos}), define 
$f(z):=C_r \sum_{j=0}^{\infty} \hat{\theta}_j z^j$ and $g(z):=\sum_{j=0}^{\infty} \tau^{j-1} z^j$. It is immediate to see that these functions satisfy 
\[
C_r z f(z) g(z)  +C_r \hat{\theta}_0=f(z) \mx{.}
\] 
Now it is sufficient to expand $f(z)$ in power of $z$, checking by induction that coefficients of the expansion (i.e. $\hat{\theta}_j$) are exactly $C_r(\tau+C_r)^{j-1}$. The argument for $\gamma_s$ is analogous.   
\endproof
Now set $\tau:=b-2 C_r$. In this way, $\tau+C_r<1$ and due to assumption (\ref{eq:minoreuno}), sequences (\ref{eq:seqces}) are monotonically decreasing. By using their expressions, the last two terms of the first of (\ref{eq:seqdue}), computed for $s+1$ can be respectively bounded in the following way 
\beq{eq:stimasecondo}
\frac{1}{s+1} \sum_{l=1}^s l h^{l-1} \theta_{s-l+1} \leq C_r
\frac{(\tau+C_r)^s}{\tau(s+1)} \sum_{l=1}^s \left( \frac{\tau}{\tau+C_r}\right)^l
\leq C_r \frac{(\tau+C_r)^s}{s+1} 
\mx{,}
\eeq
moreover, 
\beq{eq:stimaterzo}
\frac{C_r}{s+1} \sum_{l=1}^s l \beta_l h^{s-l} \leq C_r \frac{h^s}{\tau(s+1)} \sum_{l=1}^s 
\left( \frac{\tau}{h}\right)^l \leq C_r \frac{\tau^s}{3h(s+1)} \mx{,}
\eeq
as $\tau \geq 4 h$ by hypothesis. Recalling $\tau+C_r<1$, the expression for $\beta_{s+1}$ can be estimated as follows for $s \geq 1$
$$
\beta_{s+1} \leq h^s + \mu \frac{\Gamma}{s}+\frac{C_r}{(s+1)}\left(\frac{\tau^s}{3h}+1 \right) \mx{.}
$$
The inductive step is achieved once the r.h.s. of the latter is shown to be not greater than $\tau^s/(s+1)$. As $C_r/(3h) \leq 1/2 $ by hypothesis, one gets $\tau^s \geq 2 C_r + 4 \mu \Gamma + 2(s+1)h^s$, satisfied by $\tau$ chosen as above for all $s \geq 1$. This justifies the choice of $b$ as in (\ref{eq:minoreuno}).
\endproof
By using expression (\ref{eq:minoreuno}) with (\ref{eq:cr}) and (\ref{eq:gammagrande}) and recalling that $\delta<1$, the quantity appearing in (\ref{eq:minunmezzo}) can be bounded as follows 
$$
\frac{2 e \ml{F}}{d^2 \alpha \delta \sigma} + b  \leq
\frac{4 r \mu}{d^2 \alpha \delta \sigma}+
\frac{4e^2 +16(r-1)}{e d^2 \alpha \delta \sigma} \ml{F}+4h
\leq  \frac{ 4 r}{d^2 \alpha \delta \sigma} \lambda_{\mu,\ep}+4h \mx{,}
$$
in which the inequality $4e^2 +16(r-1) \leq 4 e^2 r$ for all $r \geq 1$ in which the second of (\ref{eq:ftildeeco}) have been used. It is sufficient to set $d=1/8$ and recall (\ref{eq:feh}) to show that (\ref{eq:smallness}) implies (\ref{eq:minunmezzo}) and then (\ref{eq:minoreuno}). On the other hand, by (\ref{eq:cr}), (\ref{eq:smallness}) and finally by (\ref{eq:kappa}) compared with (\ref{eq:feh})
\[
C_r \leq \frac{1}{e^2} \frac{2^{8} e r \ml{F}}{\alpha \delta \sigma}  < 
\frac{1}{e^2}   \frac{2^{8} r \lambda_{\ep,\mu}}{\alpha \delta \sigma}
 \leq 
\frac{1}{e^2}\left( \frac{1}{2}-4 h \right) \leq \frac{3}{2}h \mx{,}
\]
as required by Lemma \ref{lem:minun}. Hence (\ref{eq:psiseq}) holds and hypotheses of Lemma \ref{lem:tecnico}, statement $3$, are satisfied. A further restriction of the domain by $1/8$ completes the proof of Theorem \ref{thm:normal}.\\ 
The remainder estimate (\ref{eq:stimaresto}) is described in detail in \cite{giorgilli02}. 

\section{Nekhoroshev stability}
The Nekhoroshev type estimate contained in Theorem \ref{thm:aper} can be straightforwardly obtained combining the analytic part given by Theorem \ref{thm:normal} and the geometric part described in \cite{giorgilli02}. The main steps leading to the wanted estimate can be summarized as follows.\\
A direct consequence of Theorem  \ref{thm:normal} on a non-resonance domain, is the existence of $n-\dim \ml{M}$ \emph{approximate} first integrals. More precisely, for each independent unit vector $\lambda \in \ml{M}^{\perp}$ it is easy to check (see \cite{giorgilli02}) that the function $\Phi=T_{\chi^{(r)}} \Phi_0$ with $\Phi_0 = \langle \lambda, I \rangle$ is a first integral for the Hamiltonian in the normal form (\ref{eq:normformth}), up to the remainder $\ml{R}^{(r+1)}$. Given an initial condition $I(0) \in \ml{V}$ these first integrals determine invariant regions in phase space whose intersection is the so called \emph{plane of fast drift} $\Pi_{\ml{M}}(I(0)):=I(0)+\spn (\ml{M})$. This means that a solution starting in $I(0)$ may undergo a variation due to the resonant ``residual'' $Z$ in the normal form, which size is not controlled by the normalization order $r$, as depending only on $\ml{M}$. Hence, the action of the remainder can be interpreted as a deviating effect from this plane. However, this deviation is ``small'', as quantitatively stated in the following
\begin{cor}\label{cor:localstab}
Assume the same hypothesis of Theorem \ref{thm:normal} and consider a trajectory $(I(t),\ph(t))$ for $(\ref{eq:ham})$ such that $I(t) \in \ml{V}$ for all $t \in [a,b]$, $ab<0$. Then this solution satisfies $\dist(I(t),\Pi_{\ml{M}}(I(0))) \leq \delta/2$ for all $t \in [a,b] \cap [-t^*,t^*]$ with
\beq{eq:tempotstar}
t^*=\frac{e^2 \delta}{C_1 \ep \Delta^r} \mx{.}
\eeq  
\end{cor}
\proof Given in  \cite{giorgilli02}. \endproof 
Unfortunately, the result above holds only as long as the trajectory remains in $\ml{V}$. Roughly, the aim of the geometric part is to show that, given $I(0) \in \ml{G}$, there exists a suitable domain containing it, for which the above result can be used, then finally providing a parameters choice in a way the solution starting at $I(0)$ remains in this set for an exponentially long time. These sets, called \emph{extended block}, cover the entire action space and are shown to be non-resonance domains of the type $(\ml{M},\beta_s/2,\delta_s,N)$, where $\beta_0 <\ldots<\beta_n$ and $\delta_0 < \ldots <\delta_n <\delta/2$ are suitable sequences of parameters (see {giorgilli02}). Here hypothesis \ref{hyp:regdue} plays a key role. In this way it is possible to show that Corollary \ref{cor:localstab} acquire global validity and is enforced as follows.\\
By (\ref{eq:smallness}) and (\ref{eq:tempotstar}), define $\Delta_0:=\Delta|_{\beta=\beta_0,\delta=\delta_0}$ then $t_0^*:=t^*|_{\delta=\delta_0,\Delta=\Delta_0}$. 
\begin{prop}\label{prop:last} Assume Hypothesis \ref{hyp:regdue}.
Given $\rho>0$ suppose $\Delta_0 \leq 1/2$ for all $\delta<\rho/3$ and $K$ satisfying (\ref{eq:kappa}). Then every trajectory for $(\ref{eq:ham})$ with $I(0) \in \ml{G}$ satisfies $\dist(I(t),I(0))<\delta$ for all $|t|<t_0^*$. 
\end{prop}
By using the above mentioned values, we have
\[
\beta_0 \delta_0= K \frac{\delta^2 m^2 r}{\delta_*^2 M r^a}, \qquad \delta_*:=(n+2)! \left( \frac{4M}{m}\right)^{n+1}K^{\frac{a}{2}} \mx{,}
\]
with $a=n^2+n$. Substituting in $\Delta_0$ one gets
$$
\Delta_0=\Delta^* \frac{r^a \rho^2}{2 e \delta^2}\lambda_{\ep,\mu}+4e^{-K \frac{\sigma}{2}},\qquad 
\Delta^*:=\frac{2^{11} e \delta^* M}{\sigma K m^2 \rho^2} \mx{.}
$$
First hypothesis of Proposition \ref{prop:last} is satisfied if each addend of $\Delta_0$ is smaller than $(2e)^{-1}$ which leads to
\[
r=\left\lfloor \left(\frac{\delta^2}{\Delta^* \rho^2 \lambda_{\ep,\mu}}\right)^{\frac{1}{a}} \right\rfloor,
\qquad K= \left\lceil \Sigma  \right\rceil,\quad \Sigma:=\frac{2(1+3 \log 2)}{\sigma} \mx{,}
\] 
where $\lfloor \cdot \rfloor$ and $\lceil \cdot \rceil$ denote the rounding to the lower and to the greater integer respectively.  
The choice $\delta=(\Delta^* \lambda)^{\frac{1}{4}} \rho$ ensures that $r \geq 1$ then condition $\delta<\rho/3$ is true provided (\ref{eq:condizioneepsilon}) is satisfied. As a consequence of the above described choice for $K$, one has $K \leq 1+ \Sigma $, which satisfies (\ref{eq:kappa}) as $\sigma \leq 1 $.\\
The exponential estimate is straightforward from Proposition \ref{prop:last}, by replacing the already determined expression for $r$ and $\Delta_0 \leq 1/e$ in $t_0^*$. This completes the proof of Theorem \ref{thm:aper}. 

\subsection*{Acknowledgements} We would like to acknowledge useful e-mail exchanges with Prof. Antonio Giorgilli.

\bibliographystyle{alpha}
\bibliography{rev_aper_Nek}

\def\cprime{$'$}
\begin{thebibliography}{BGG85}

\bibitem[Arn63]{arn1}
V.~I. Arnold.
\newblock Proof of {A}. {N}. {K}olmogorov's theorem on the preservation of
  quasiperodic motions under small perturbations of the {H}amiltonian.
\newblock {\em Russ. Math. Surveys}, 18(5):9--36, 1963.

\bibitem[BGG85]{bgg85}
G.~Benettin, L.~Galgani, and A.~Giorgilli.
\newblock A proof of {N}ekhoroshev's theorem for the stability times in nearly
  integrable {H}amiltonian systems.
\newblock {\em Celestial Mech.}, 37:1--25, 1985.

\bibitem[Bou13]{boun13}
A.~Bounemoura.
\newblock Effective stability for slow time-dependent near-integrable
  {H}amiltonians and application.
\newblock {\em C. R. Math. Acad. Sci. Paris}, 351(17-18):673--676, 2013.

\bibitem[Gal86]{galla}
G.~Gallavotti.
\newblock Quasi-integrable mechanical systems.
\newblock In {\em Ph\'enom\`enes critiques, syst\`emes al\'eatoires, th\'eories
  de jauge, {P}art {I}, {II} ({L}es {H}ouches, 1984)}, pages 539--624.
  North-Holland, Amsterdam, 1986.

\bibitem[Gio02]{giorgilli02}
A.~Giorgilli.
\newblock Notes on exponential stability of {H}amiltonian systems.
\newblock In {\em Dynamical Systems. Part I. Hamiltonian Systems and Celestial
  Mechanics}, Pisa, 2002. Centro di Recerca Matematica Ennio De Giorgi, Scuola
  Normale Superiore.

\bibitem[Gli64]{gli}
J.~Glimm.
\newblock Formal stability of {H}amiltonian systems.
\newblock {\em Comm. Pure Appl. Math.}, 17:509--526, 1964.

\bibitem[GZ92]{gz92}
A.~Giorgilli and E.~Zehnder.
\newblock Exponential stability for time dependent potentials.
\newblock {\em Z. angew. Math. Phys. (ZAMP)}, 43:827--855, 1992.

\bibitem[Kol54]{kolm}
A.~N. Kolmogorov.
\newblock On conservation of conditionally periodic motions under small
  perturbations of the {H}amiltonian.
\newblock {\em Dokl. Akad. Nauk. USSR}, 98(4):527--530, 1954.

\bibitem[Lit59a]{littlewood59b}
J.~E. Littlewood.
\newblock The {L}agrange configuration in celestial mechanics.
\newblock {\em Proc. London Math. Soc.}, 9(3):525--543, 1959.

\bibitem[Lit59b]{littlewood59a}
J.~E. Littlewood.
\newblock On the equilateral configuration of the restricted problem of three
  bodies.
\newblock {\em Proc. London Math. Soc.}, 9(3):343--372, 1959.

\bibitem[Loc92]{lo}
P.~Lochak.
\newblock Canonical perturbation theory: an approach based on joint
  approximations.
\newblock {\em Uspekhi Mat. Nauk}, 47(6(288)):59--140, 1992.

\bibitem[MG95]{morbidelli}
A.~Morbidelli and A.~Giorgilli.
\newblock Superexponential stability of {KAM} tori.
\newblock {\em J. Statist. Phys.}, 78(5-6):1607--1617, 1995.

\bibitem[Mos55]{moser55}
J.~Moser.
\newblock Stabilit\"atsverhalten kanonischer {D}ifferentialgleichungssysteme.
\newblock {\em Nachrichten der Akademie der Wissenschaften in G\"ottingen. II.
  Mathematisch-Physikalische Klasse}, 1955:87--120, 1955.

\bibitem[Nek77]{nekhoroshev77}
N.~N. Nekhoroshev.
\newblock An exponential estimate on the time of stabilty of nearly-integrable
  {H}amiltonian systems.
\newblock {\em Russ. Math. Surveys}, 32:1--65, 1977.

\bibitem[Poi92]{poi}
H.~Poincar{\'e}.
\newblock {\em Les m\'ethodes nouvelles de la m\'ecanique c\'eleste.}
\newblock Gauthier-Villars, Paris, 1892.

\bibitem[P{\"o}s93]{poes}
J.~P{\"o}schel.
\newblock Nekhoroshev estimates for quasi-convex {H}amiltonian systems.
\newblock {\em Math. Z.}, 213(2):187--216, 1993.

\bibitem[P{\"o}s01]{poes2}
J{\"u}rgen P{\"o}schel.
\newblock A lecture on the classical {KAM} theorem.
\newblock In {\em Smooth ergodic theory and its applications ({S}eattle, {WA},
  1999)}, volume~69 of {\em Proc. Sympos. Pure Math.}, pages 707--732. Amer.
  Math. Soc., Providence, RI, 2001.

\bibitem[PW94]{perry}
A.~D. Perry and S.~Wiggins.
\newblock K{AM} tori are very sticky: rigorous lower bounds on the time to move
  away from an invariant {L}agrangian torus with linear flow.
\newblock {\em Phys. D}, 71(1-2):102--121, 1994.

\bibitem[WM13]{wiggins}
S.~Wiggins and A.~Mancho.
\newblock Barriers to transport in aperiodically time-dependent two-dimensional
  velocity fields: Nekhoroshev's theorem and ``nearly invariant'' tori.
\newblock 2013.

\end{thebibliography}

\end{document}